\documentclass[12pt, 14paper,reqno]{amsart}
\usepackage{amsmath,amsfonts,amssymb}
\usepackage{mathrsfs}
\usepackage{longtable}
\usepackage[breaklinks]{hyperref}
\usepackage{graphicx}
\makeatletter
\@namedef{subjclassname@2020}{%
  \textup{2020} Mathematics Subject Classification}
\makeatother
\theoremstyle{plain}
\newtheorem{thm}{Theorem}[section]
\newtheorem{lem}{Lemma}[section]

\newtheorem{thma}{Theorem}

\theoremstyle{proof}
\numberwithin{equation}{section}

\begin{document}
\title[Fields with small class group in the family $\mathbb{Q}(\sqrt{9m^2+2m})$]{ 
Fields with small class group in the family $\mathbb{Q}(\sqrt{9m^2+2m})$}
\author{Kalyan Chakraborty and Azizul Hoque}
%\address{K. Banerjee @Department of Mathematics, SRM University AP, Mangalagiri-Mandal, Amaravati-522502, Andhra Pradesh, India.} \email{kalyan.ba@srmap.edu.in}
\address{K. Chakraborty @Department of Mathematics, SRM University AP, Mangalagiri-Mandal, Amaravati-522502, Andhra Pradesh, India.} \email{kalyan.c@srmap.edu.in}
\address{A. Hoque @Department of Mathematics, Faculty of Science, Rangapara College, Rangapara, Sonitpur-784505, Assam, India.}
\email{ahoque.ms@gmail.com}

\maketitle

\begin{abstract} Very recently, Issa and Darrag [Arch. Math. (Basel) {\bf 123} (2024), no. 4, 379--383] determined partial Dedekind zeta values for certain ideal classes in  the real quadratic fields of the form  $\mathbb{Q}(\sqrt{9m^2+2m})$, where $9m^2+2m$ is square-free and $m\equiv 2\pmod 3$ is an odd positive integer. We use these partial Dedekind zeta values to investigate the small class numbers of such fields. More precisely, we prove that the class numbers of the fields in the above mentioned family are at least $4$. Further, we provide a sufficient condition permitting to specify the structure of the class groups of order $4$ in this family of fields.
\end{abstract}

\section{Introduction}
It was conjectured by Gauss' that there are infinitely many real quadratic fields with class number $1$. This conjecture is still open. There are some interesting results (cf. \cite{BI03, BL16, BK1, KKST19, LO88, MO87}) along this conjecture for certain families of real quadratic fields and most of these families have only a few fields with class number $1$. This problem is resolved for a well known family; so-called  Richaud-Degert (R-D) type fields of real quadratic fields since their fundamental units are explicitly known (see, \cite{DE58}). These fields are of the form $\mathbb{Q}(\sqrt{d})$, where $d=n^2+r$ is square-free with $n, r\in \mathbb{Z}, -r<n\leq r$ and $r\mid 4n$.   
There are some interesting results for small class numbers of R-D type fields. For example, Byeon and Kim \cite{BK2} determined some criteria for the class number $2$ of these fields. Chakraborty, Hoque and Mishra \cite{CHM19} extended the results of Byeon and Kim \cite{BK2}, and they determined some criteria for the class number up to $3$ when $r=\pm 1, \pm 4$. Moreover, they proved in \cite{CHM21} that the class number of this family of fields with $r=1$ can be made as large as possible by increasing the number of prime divisors of $n$. Biro and Lapkova \cite{BL16} determined the fields with class number $1$ in the family $\mathbb{Q}(\sqrt{9m^2+4m})$. Recently, Mahapatra, Pandey and Ram \cite{MPR23} extended the work of Biro and Lapkova \cite{BL16}, and proved that there is only one field with class number $1$ and one field with class number $2$ in the family of fields $\mathbb{Q}(\sqrt{9m^2+4m})$ when $m\equiv 1 \pmod 3$. 

Very recently, Issa and Darrag \cite{ID24} determined partial Dedekind zeta values for certain ideal classes in the real quadratic fields of the form  $k_m:=\mathbb{Q}(\sqrt{9m^2+2m})$, where $9m^2+2m$ is square-free and $m\equiv 2\pmod 3$ is an odd positive integer. We use these partial Dedekind zeta values to investigate such fields having small class numbers. Gauss' genus theory demonstrates that the class number, $h(d)$ of $\mathbb{Q}(\sqrt{d})$ is odd if and only if $d = p, 2p_1$, or $p_1 p_2$ where $p$ is prime and $p_1 \equiv p_2 \equiv 3 \pmod 4$ are primes. Since $m\equiv 2\pmod 3$ is odd, so that $h(9m^2+2m)$ is even and thus $h(9m^2+2m)>1$. The same result follows from the paper of Byeon and Kim \cite[Corollary 3.8]{BK1}. Moreover, Corollary 3.5 in \cite{BK2} gave certain criteria for $h(9m^2+2m)$ to be $2$. In this paper, we prove that $h(9m^2+2m)\geq 4$. Further, we provide a sufficient condition permitting to specify the structure of the class groups of order $4$ in this family of fields. The precise result is the following:
\begin{thm}\label{thm1}
Let $m\equiv 2\pmod 3$ be an odd positive integer. If $9m^2+2m$ is square-free, then $h(9m^2+2m)\geq 4$. Moreover, if $h(9m^2+2m)=4$, then the following hold:
\begin{itemize}
\item[(i)] $\text{Cl}(k_m)\cong \mathbb{Z}/2\mathbb{Z}\times \mathbb{Z}/2\mathbb{Z}$,
\item[(ii)] for a prime divisor $p\geq 5$ of $m$,  we have %$\sum\limits_{\substack{ |t|<\sqrt{9m^2+2m}\\ t^2\equiv 4(9m^2+2m)\pmod 4}} \sigma(\frac{4(9m^2+2m)-t^2}{4}),$
$$\sum\limits_{|\ell|<\sqrt{9m^2+2m}} \sigma(9m^2+2m-\ell^2)=\frac{1}{3}(147m^3+49m^2+308m+46)+\frac{m^2}{p^2}(36m+12) + 9mp^2+10m+9p^2,$$
where $\sigma(*)$ denotes the sum of divisors of $*$. 
\end{itemize}
\end{thm}
We first use partial Dedekind zeta values at $-1$ of certain ideal classes in the family of fields $\mathbb{Q}(\sqrt{9m^2+2m})$ to show that the class number of such fields are at least $4$. To classify the class groups of oder $4$, we compare these values with complete Dedekind zeta value at $-1$ computed with the help of classical Siegel's formula as described by D. B. Zagier \cite{ZA76}. 
\section{Partial Dedekind zeta values}
We first mention a classical result of Lang \cite{LAN} for calculating the special values of Dedekind zeta function,  $\zeta_k$ when $k$ is real quadratic field. Given an ideal class $\mathfrak{A}$ of $k$,  let $\mathfrak{a}$ be an integral ideal in $\mathfrak{A}^{-1}$ with  an integral basis $\{r_{1},r_{2}\}$, i.e.  $\mathfrak{a}= \mathbb{Z}r_1 \oplus \mathbb{Z}r_2$, where $r_1,r_2 \in \mathfrak{a}$. We set 
$$\delta(\mathfrak{a}):= r_1r_2'-r_1'r_2,$$ 
where $r_1'$ and  $r_2'$ are the conjugates of $r_1$ and $r_2$ 
respectively.

Let $\varepsilon$ be the fundamental unit of $k$. Then $           
\{\varepsilon r_1, \varepsilon r_2\}$ is also an integral basis of 
$\mathfrak{a}$, and thus we can find a matrix 
$M=
\begin{bmatrix}
a&b\\
c&d
\end{bmatrix}
$
with integer entries satisfying the following: 
$$
\varepsilon\begin{bmatrix}
r_1\\r_2
\end{bmatrix}
=M\begin{bmatrix}
r_1\\
r_2
\end{bmatrix}.
$$
The classical result of Lang which helps us to compute partial Dedekind zeta  value for $\mathfrak{A}$ at $-1$ is the following:

\begin{thma}[{\cite[p. 159]{LAN}}]\label{thmlang}
By keeping the above notations, we have 
\begin{align*}
\zeta_k(-1, \mathfrak{A})&={\displaystyle\frac{\textsl{sgn }
\delta(\mathfrak{a})~r_2r_2'}{360N(\mathfrak{a})c^3}}\big\{(a
+d)^3-6(a+d)N(\varepsilon)-240c^3(\textsl{sgn } c)\\
&\times S^3(a,c)+180ac^3(\textsl{sgn } 
c)S^2(a,c)-240c^3(\textsl{sgn } c)S^3(d,c)\\
& +180dc^3(\textsl{sgn } c)S^2(d,c) \big\},
\end{align*}
where $N(\mathfrak{a})$ represents the norm of $\mathfrak{a}$ and 
$S^i(-,-)$ denotes the generalized Dedekind sum as defined in 
\cite{AP50}.
\end{thma} 

We need to determine the values of $a,b,c,d$ and generalized 
Dedekind sums in order to apply Theorem \ref{thmlang}. The following result (see, \cite[p. 143, Eq. 2.15]{LAN}) helps us to determine these values.

\begin{lem}\label{lem1} The matrix $M$ is given by 
$$
\begin{bmatrix}
Tr\left(\frac{r_1 r_2'\varepsilon}{\delta(\mathfrak{a})}\right)& Tr
\left(\frac{r_1 r_1'\varepsilon'}{\delta(\mathfrak{a})}\right) 
\vspace*{2mm} \\ 
 Tr\left(\frac{r_2 r_2'\varepsilon}{\delta(\mathfrak{a})}\right) &
 Tr\left(\frac{r_1 r_2'\varepsilon'}{\delta(\mathfrak{a})}\right)
\end{bmatrix}.
$$
Moreover, $\det(M)=N(\varepsilon)$ and $bc\ne 0$.
\end{lem} 
In order to use Theorem \ref{thmlang}, one needs to calculate generalized Dedekind sums. In this particular fields, we need the following sums \cite[p. 155; Eq. (4.3)-(4.4)]{LAN}.
\begin{lem}\label{lem2} For any positive integer $m$, we have
\begin{itemize}
\item[(i)] $S^3(\pm 1, m)={\displaystyle\pm\frac{-m^4+5m^2-4}{120m^3}}.$
\vspace*{2mm}
\item[(ii)] $S^2(\pm 1, m)={\displaystyle\frac{m^4+10m^2-6}{180m^3}}.$
\end{itemize}
\end{lem}
%We determine Dedekind zeta values attached to $k_n$ in two ways 
%using Theorem \ref{thm2.1} and Theorem \ref{thm2.2} and then 
%compare these values and finally use elementary group theoretic 
%arguments to establish our results.  

We are now in a position to discuss the partial Dedekind zeta values at $-1$ for certain ideal classes in the field $k_m:=\mathbb{Q}(\sqrt{9m^2+2m})$. Byeon and Kim \cite[Theorem 2.3]{BK1} used Theorem \ref{thmlang} together with Lemmas \ref{lem1} and \ref{lem2} to compute the partial Dedekind zeta value for the principal ideal class $\mathcal{P}$ of R-D type real quadratic fields. By this theorem, we have 
\begin{equation}\label{eqP}
\zeta_{k_m}(-1,\mathcal{P})=\frac{108m^3+36m^2+57m+9}{180},
\end{equation}
where $\mathcal{P}$  denotes the principal ideal class in $k_m$. 

It is easy to see that $2$ ramifies in the field $k_m$ and thus  
$$(2)=(2,1+\sqrt{9m^2+2m})^2.
$$
Let $\mathcal{A}$ be the ideal class containing  $(2,1+\sqrt{9m^2+2m})$. Then by \cite[Theorem 2.5 (I)(i)]{BK1}, one gets
\begin{equation}\label{eqA}
\zeta_{k_m}(-1,\mathcal{A})=\frac{27m^3+9m^2+138m+36}{180}.
\end{equation}

Since $3$ splits in $k_m$ as 
$$(3)=(3,1+\sqrt{9m^2+2m})(3,1-\sqrt{9m^2+2m}),$$ 
if $\mathcal{B}$ is the ideal class containing  $(3,1+\sqrt{9m^2+2m})$  then $(3,1-\sqrt{9m^2+2m})\in \mathcal{B}^{-1}$. Thus  $\zeta_{k_m}(-1,\mathcal{B})=\zeta_{k_m}(-1,\mathcal{B}^{-1})$. These values are evaluated by Issa and Darrag \cite[Theorem 1]{ID24}, which are 
\begin{equation}\label{eqB}
\zeta_{k_m}(-1,\mathcal{B})=\zeta_{k_m}(-1,\mathcal{B}^{-1})=\frac{12m^3+4m^2+113m+1}{180}. 
\end{equation}
Further, for any prime divisor $p>3$ of $m$,  $p$ remifies in $k_m$ and $(p)=(p, \sqrt{9m^2+2m})^2$. Let $\mathcal{C}$ be the ideal class containing  $(p,\sqrt{9m^2+2m})$. The value of $\zeta_{k_m}(-1, \mathcal{C})$ is obtained by Issa and Darrag \cite[Theorem 2]{ID24} as
\begin{equation}\label{eqC}
\zeta_{k_m}(-1, \mathcal{C})=\frac{108m^3+36m^2+27p^4m+30p^2m+9p^4}{180p^2}.
\end{equation}

\section{Proof of Theorem \ref{thm1}}
We first compare the values of $\zeta_{k_m}(-1,\mathcal{P})$ and $\zeta_{k_m}(-1,\mathcal{A})$, that is,  $\zeta_{k_m}(-1,\mathcal{P})=\zeta_{k_m}(-1,\mathcal{A})$ gives us $m=1$, which is not possible. Thus, $\mathcal{A}$ is non-principal ideal class in $\mathcal{O}_{k_m}$. Moreover, $(2,1+\sqrt{9m^2+2m})^2=(2)$ and $(2,1+\sqrt{9m^2+2m})\in \mathcal{A}$ together imply $|\mathcal{A}|=2$, and hence $h(9m^2+2m)$ is even. Similarly, we can verify that $\mathcal{B}$ and $\mathcal{C}$ are also non-principal ideals in $\mathcal{O}_{k_m}$. 

We now pairwise compare the values of $\zeta_{k_m}(-1,\mathcal{A}), \zeta_{k_m}(-1,\mathcal{B})$ and $\zeta_{k_m}(-1,\mathcal{C})$ to see the distinctness. If $\zeta_{k_m}(-1,\mathcal{A})=\zeta_{k_m}(-1,\mathcal{B})$, then $3m^3+m^2+5m+7=0$, which is not possible for any positive integer $m$. Therefore, $\mathcal{A}$ and  $\mathcal{B}$ are distinct. Along the same lines, we can verify that $\mathcal{C}$ is distinct from  both $\mathcal{A}$ and  $\mathcal{B}$. This confirms that there are at least $4$ distinct ideal classes in $\mathcal{O}_{k_m}$, and therefore $h(9m^2+2m)\geq 4$.  

We now assume that $h(9m^2+2m)= 4$. Then $\mathcal{P}, \mathcal{A}, \mathcal{B}$ and $\mathcal{C}$ are the only ideal classes in $\mathcal{O}_{k_m}$. However, we have seen that $\mathcal{A}$ is of order $2$. Since $(p, \sqrt{9m^2+2m})^2=(p)$ and $(p, \sqrt{9m^2+2m})\in \mathcal{C}$, so that $\mathcal{C}$ is also of order $2$. We see from \eqref{eqB} that 
$\zeta_{k_m}(-1,\mathcal{B})=\zeta_{k_m}(-1,\mathcal{B}^{-1})$. On the other hand, $\zeta_{k_m}(-1,\mathcal{B}^{-1})\ne \zeta_{k_m}(-1,\mathcal{A})\ne \zeta_{k_m}(-1,\mathcal{C})$. Thus, $\mathcal{B}=\mathcal{B}^{-1}$ and hence the order of $\mathcal{B}$ is also $2$. Therefore $Cl(k_m)\cong\mathbb{Z}/2\mathbb{Z}\times \mathbb{Z}/2\mathbb{Z}$. 

We now recall a particular case of a celebrated result due to Siegel \cite{SI69} about the computation of Dedekind zeta value, $\zeta_k(1-2n)$, where $n$ is a positive integer and $k$ is a real quadratic field. 
\begin{thma}[{\cite[p.76]{ZA76}}]\label{thmazeta1}
Let $k$ be a real quadratic field with discriminant $D$. Then 
$$\zeta_k(-1)=\frac{1}{60}\sum_{\substack{|b|<\sqrt{D}\\ b^2\equiv D\pmod 4}}  \sigma\left(\frac{D-b^2}{4}\right),$$
where $\sigma(*)$ denotes the sum of divisors of $*$. 
\end{thma}

Since $\mathcal{P}, \mathcal{A}, \mathcal{B}$ and $\mathcal{P}$ are the only ideal classes in $\mathcal{O}_{k_m}$, so that
$$\zeta_{k_m}(-1)=\zeta_{k_m}(-1, \mathcal{P})+\zeta_{k_m}(-1, \mathcal{A})+\zeta_{k_m}(-1, \mathcal{B})+\zeta_{k_m}(-1, \mathcal{C}).$$
Therefore by Theorem \ref{thmazeta1}, we get
$$\frac{1}{60}\sum_{\substack{|b|<\sqrt{D}\\ b^2\equiv D\pmod 4}}  \sigma\left(\frac{D-b^2}{4}\right)=\zeta_{k_m}(-1, \mathcal{P})+\zeta_{k_m}(-1, \mathcal{A})+\zeta_{k_m}(-1, \mathcal{B})+\zeta_{k_m}(-1, \mathcal{C}).$$
Since $D=4(9m^2+2m)$, $b$ must be even, and thus we write $b=2\ell$. Then the above equation becomes,
$$\frac{1}{60}\sum_{\substack{|\ell|<\sqrt{9m^2+2m}}}  \sigma(9m^2+2m-\ell^2)=\zeta_{k_m}(-1, \mathcal{P})+\zeta_{k_m}(-1, \mathcal{A})+\zeta_{k_m}(-1, \mathcal{B})+\zeta_{k_m}(-1, \mathcal{C}).$$

Utilizing \eqref{eqP}--\eqref{eqC}, we get
$$\sum\limits_{|\ell|<\sqrt{9m^2+2m}} \sigma(9m^2+2m-\ell^2)=\frac{1}{3}(147m^3+49m^2+308m+46)+\frac{m^2}{p^2}(36m+12) + 9mp^2+10m+9p^2.$$
This completes the proof of Theorem \ref{thm1}. 
%We end this article with the remark below. 
%\begin{remark}
%In place of $9m^2+2m$, if we consider $p^2m^2+2m$ with $m\equiv 2\pmod p$ for any prime $p\geq 3$, then the similar result follows. 
%\end{remark}

\section*{Acknowledgements} 
\noindent 
%The authors would like to express their gratitude to Professor Claude Levesque for carefully reading this manuscript and for his useful comments. The second author is grateful to Professor Srinivas Kotyada for stimulating environment at The Institute of Mathematical Sciences, Chennai during his visiting period. The authors are thankful to the anonymous referees for their valuable  comments and suggestions which have helped improving the presentation immensely. 
This work is supported by ANRF-SERB MATRICS grant (No. MTR/2021/000762) and ANRF-SERB CRG grant (No. CRG/2023/007323), Govt. of India. %%%%%%%%%%%%%%%%%%%%%%%%%%%%%%%%%%%%%%%%%%%%%%%%%%%      

\end{document}